\documentclass[11pt,reqno]{amsart}

\usepackage{amsmath}
\usepackage{amssymb}
\usepackage{graphicx,color}
\usepackage{algorithm}
\usepackage{algpseudocode}

\usepackage[utf8]{inputenc}
\usepackage{listings}
\usepackage{xcolor}

\usepackage[utf8]{inputenc}
\usepackage{listings}
\usepackage{xcolor}

\usepackage[utf8]{inputenc}
\usepackage{listings}
\usepackage{xcolor}

\usepackage[margin=2.4cm]{geometry}

\newtheorem{theorem}{Theorem}[section]
\newtheorem{corollary}[theorem]{Corollary}

\newtheorem{proposition}[theorem]{Proposition}
\theoremstyle{definition}

\DeclareMathOperator*{\id}{id}

\numberwithin{equation}{section}

\title[Distribution steering in discrete time]{A first-order condition  for discrete-time  distribution steering}

\author[A. Dom\'inguez Corella]{Alberto Dom\'inguez Corella}

\address[A. Dom\'inguez Corella]{Institut f\"{u}r Stochastik und Wirtschaftsmathematik, Variational Analysis, Dynamics and Operations Research Unit E105-04, Technische Universit\"at Wien, Wiedner Hauptstra\ss{e} 8, 1040 Vienna, Austria}
\address{Institut f\"ur Mathematik und Wissenschaftliches Rechnen,
Universit\"at Graz, Heinrichstrasse 36, 8010 Graz, Austria}
\email{{\tt alberto.of.sonora@gmail.com}}

\author[D. Gonz\'alez-S\'anchez]{David González-Sánchez}
\address[D. Gonz\'alez-S\'anchez]{SECIHTI--CINVESTAV-IPN, Mathematics Department, Mexico City 07000, Mexico}
\email{\tt dgonzalezsa@secihti.mx}

\keywords{distribution steering, discrete-time, optimal control, Pontryagin}

\subjclass[2010]{49J55, 49K45,49K40,90C40}

\begin{document}

	\begin{abstract}
		We study a class of distribution-steering problems from a variational point of view. Under some differentiability assumptions, we derive necessary conditions for optimal Markov policies in the spirit of the Lagrange multiplier approach. We also provide a heuristic gradient-based method derived from the variational principle.
	\end{abstract}
	
	\maketitle
	
	
	\section{Introduction}
	Distribution steering problems are a class of control problems in which the system's state is represented by a probability distribution, and the objective is to steer the initial distribution toward a specified target distribution. This class of problems, in discrete time, have recently attracted significant attention due to its theoretical foundations and wide range of applications \cite{Bakolas_2016, Bakolas_2017, Bakolas_SCL_2018, Bakolas_Automatica_2018, BalciBakolas_2021,deBaydinEtAl,  ItoKashima_2025, TerpinEtAl}.
	\smallbreak 
	In this paper, we give a variational principle for a stochastic optimal control problem in discrete time associated with distribution steering. The variational principle we provide is a type of stochastic Pontryagin principle. It can also be seen as a necessary condition involving Lagrange multipliers.
	\smallbreak 
	We deal with a discrete-time dynamical system, with dynamics $f_t:\mathbb R^n\times \mathbb R^m\to\mathbb R^n$, that is controlled by Markov policies $\varphi_t:\mathbb R^n\to \mathbb R^m$. Given an initial random variable $x_0$, with values in $\mathbb R^n$, there is a finite sequence $x_0,\ldots,x_T$ of states, where 
	\[ x_{t+1}=f_t\big(x_t,\varphi_t(x_t)\big),\qquad t=0,1,\ldots,T-1.
	\]
	The \textit{distribution steering problem} we consider consists of finding a Markov policy $\varphi=\{\varphi_t\}_{t=0}^{T-1}$ that minimizes  the loss function
	\begin{equation}\label{FirstCostFunct}
		\frac{1}{2}\hspace*{0.02cm}\mathbb E|x_T-\mathfrak t(x_0)|^2,    
	\end{equation}
	where $\mathfrak t(x_0)$ represents a random variable with the target distribution (which is implicitly assumed to be the pushforward of the initial distribution under a known  transport map). Under differentiable dynamics and other suitable hypotheses, we derive necessary conditions for optimal Markov policies,  extending classical variational techniques to the space of probability distributions. We also propose a  gradient-based method to numerically approximate control policies that steer (in the sense of minimizing \eqref{FirstCostFunct}) the  initial distribution to the target one.
	\smallbreak 
	Some similar problems have been studied in recent years through optimal transport and optimal control. Some of the first works in discrete time include \cite{Bakolas_2016, Bakolas_2017, Bakolas_SCL_2018, Bakolas_Automatica_2018, BalciBakolas_2021}, where stochastic linear systems are considered. More recently, in \cite{deBaydinEtAl}, the authors make a connection between optimal transport and
	stochastic control problems representing distribution steering; they  find closed-form solutions for a transport map representing optimal policies. In \cite{TerpinEtAl}, the authors study discrete-time dynamic programming problems in probability spaces—which particularly model distribution steering—by means of an optimal transport approach. We also mention \cite{ItoKashima_2025}, where distribution steering is studied together with entropy maximization.
	\smallbreak 
	In contrast with previous contributions in discrete time, we follow a variational approach based on the closed-loop maximum principle in \cite[Theorem 23]{CorellaHL_2019}, which has its roots in the open-loop maximum principle from \cite[Theorem 2.2]{Aseev_2017} (see also \cite{KrastanovStefanov_2024}). These principles have the advantage of admitting a closed-form expression for the adjoint variable and have been shown effective in finding Nash equilibria in dynamic games; see, e.g., \cite{AlNassirLafta_2021,KeskinSaglam_2022}.
	\smallbreak 
	We also mention some contributions in continuous time. Constructive results are available for linear systems with Gaussian distributions, where affine control laws restrict the evolution to the space of Gaussians, allowing closed-form solutions; see, e.g., \cite{ChenGeorgiouPavon_JOTA_2016,ChenGeorgiouPavon_TAC_PartI_2016, ChenGeorgiouPavon_TAC_PartII_2016, ChenGeorgiouPavon_TAC_PartIII_2018}. These works rely on classical tools from optimal control. In \cite{Bonnet_2019, BonnetFrankowska,BonnetRossi_2019}, the authors develop a Pontryagin maximum principle for optimal control problems in the Wasserstein space, combining techniques from optimal control with the differential structure of Wasserstein geometry.
	\smallbreak 
	The paper is organized as follows. This introduction ends with the notation used throughout the manuscript. In the next section, we present the precise formulation of the problem and state our results. We also provide a heuristic gradient-based method and illustrate it with two examples. The proofs of our results are given in Section 3. 
	\smallbreak 
	{\sc Notation.} The set $\mathbb R^d$ denotes the Euclidean space equipped with its usual norm $|\cdot|$ and canonical inner product $\langle\cdot,\cdot\rangle$. The elements of $\mathbb R^d$ are regarded as column vectors. 
	
	Given a differentiable function $f:\mathbb R^n\times\mathbb R^m\to \mathbb R^n$, the partial gradients $\frac{\partial f}{\partial x}(x,u)\in M_{n\times n}(\mathbb R)$ and $\frac{\partial f}{\partial u}(x,u)\in M_{n\times m}(\mathbb R)$  at a point $(x,u)$ are given by 
	\begin{align*}
		\Big[	\frac{\partial f}{\partial x}(x,u) \Big]_{ij}:= \frac{\partial f_i}{\partial x_{j}}(x,u)\quad \text{and}\quad 	\Big[	\frac{\partial f}{\partial u}(x,u) \Big]_{ij}:= \frac{\partial f_i}{\partial u_{j}}(x,u).
	\end{align*}
	
	The normal cone to the convex set $U\subset \mathbb R^m$ at $\hat{x}\in U$ is
	\[ N_U(\hat{x}):=\{ y\in \mathbb R^m \mid \langle y, x-\hat{x} \rangle \leq 0 \quad\forall x\in U \} .  \]
    Given two  sets $A,B\subset\mathbb R^m$, their sum is defined as
    \begin{align*}
        A+B = \{a+b:\,\,\, a\in A\,\,\, \text{and}\,\,\,b\in B\}.
    \end{align*}
	Let $M^\top$ denote the transpose of the matrix $M$. We also denote the product of matrices 
	\[ \prod_{s=\tau}^{t}M_{s} =\begin{cases}
		M_\tau\cdots M_t & \text{if } t\geq \tau,\\
		I & \text{if } t<\tau.
	\end{cases}\]
	
	The set of all Borel maps $\mathfrak t: \mathbb R^n \to \mathbb R^n$ such that $$\int_{\mathbb R^n}|\mathfrak t(x)|^2 \mu(dx)<\infty$$ is denoted $L^2(\mu)^n$. Finally,  $\mathcal P_2(\mathbb R^n)$ consists of all the Borel measures on $\mathbb R^n$ with finite second moment, i.e., such that the identity mapping $\id:\mathbb R^d\to\mathbb R^d$ belongs to $L^2(\mu)$.

	\section{Problem formulation and results}
	
	\subsection{The model}
	Let $T\in\mathbb N$ denote a time horizon.  The states evolve in the Euclidean space $\mathbb R^n$ and the controls in a nonempty Borel set $U\subset \mathbb R^m$. 
	\smallbreak 
	
	A  \textit{Markov policy} is 
	a  sequence $\{\varphi_t\}_{t=0}^{T-1}$ of Borel  functions $\varphi_t:\mathbb R^n\to \mathbb R^m$ satisfying
	\begin{align*}
		\varphi_t(x)\in U\quad\forall x\in \mathbb R^n.
	\end{align*}
	The set of all Markov policies is denoted by $\Phi$. 
	The evolution of the system is represented by a sequence $\{f_t\}_{t=0}^{T-1}$ of  functions $f_t:\mathbb R^n\times \mathbb R^m\to \mathbb R^n$.
	
	\smallbreak 
	The initial state $x_0$ is a random variable with distribution $\mu\in\mathcal P_2(\mathbb R^n)$.
	Given a Markov policy $\varphi=\{\varphi_t\}_{t=0}^{T-1}$, the associated sequence of states $\{x_{t}^\varphi\}_{t=0}^{T}$ is given by $x_0^\varphi=x_0$ and 
	\begin{align}\label{system}
		x^\varphi_{t+1}=f_t\big(x_t^\varphi,\varphi_t(x_t^\varphi)\big).
	\end{align}
	For a given a map  $\mathfrak t\in L^2(\mu)^n$,  the  objective functional $\mathcal J:\Phi\to\mathbb R\cup\{+\infty\}$  is given by 
	\begin{align}\label{cost}
		\mathcal J(\varphi):=\frac{1}{2}\hspace*{0.02cm}\mathbb E|x_T^{\varphi}-\mathfrak t(x_0)|^2 .
	\end{align}
	We consider the \textit{distribution steering problem}
	\begin{align*}
		(\mathcal P)\quad \min_{\varphi\in\Phi} \mathcal J(\varphi).
	\end{align*}
	The idea behind problem $(\mathcal P)$ is to find a Markov policy such that system (\ref{system}) approximately steers $\mu$ to $\mathfrak{t}_{\#}\mu$. 
	A policy $\hat \varphi\in\Phi$ is said to be optimal for problem $(\mathcal P)$ if $\mathcal J(\hat\varphi)<+\infty$ and 
	\begin{align*}
		\mathcal J(\hat\varphi) \le \mathcal J(\varphi)\quad \forall\varphi\in \Phi. 
	\end{align*}
	
	\begin{proposition}\label{Propp}
		For each $t\in\{0,\dots, T-1\}$, suppose that 
		\begin{itemize}
			\item[(i)]  for all $x\in\mathbb R^n$ and $r>0$, the set $\big\{u\in U:\, |f_t(x,u)|\le r\big\}$ is compact; 
			
			\item[(ii)] the function $f_t:\mathbb R^n\times\mathbb R^m\to \mathbb R^n$ is continuous and there exists $m_t>0$ such that for all $x\in\mathbb R^d$ there exists $u_x\in U$ satisfying
			\begin{align*}
				|f_t(x,u_x)|\le m_t\big(1+|x|\big).
			\end{align*}
		\end{itemize}
		Then, there exists $\varphi\in\Phi$ such that $\mathcal J(\varphi)<\infty$.  
	\end{proposition}

	\subsection{A variational principle}
	We now state the announced variational principle. This can also be seen as a Pontryagin principle for discrete-time stochastic tracking problems.
	\begin{theorem}\label{Pontryagin}
		Let the following statements hold. 
		\begin{itemize}
			\item[(i)] The control  set  $U\subset\mathbb R^m$  is  convex. 
			
			\item[(ii)]  For each $t\in\{0,\dots,T-1\}$, the function $f_t:\mathbb R^n\times\mathbb R^m\to \mathbb R^n$ is differentiable  with uniformly continuous partial derivatives and there exists 
			$m_t>0$ such that
			\begin{align*}
				|f_t(x,u)|\le m_t\big(1 + |x| + |u|\big)\quad \forall (x,u)\in \mathbb R^n\times\mathbb R^m.
			\end{align*}
		\end{itemize}
		Let $ \hat\varphi\in\Phi$ be such that each $\hat\varphi_t:\mathbb R^n\to\mathbb R^m$ is differentiable with uniformly continuous partial derivatives.  Consider the sequence $\left\lbrace {p}_t\right\rbrace_{t=1}^{T}$ of  random variables given by
		\begin{align}\label{adjoint}
			p_{t}:=	\prod_{s=t}^{T-1}\Big(\frac{\partial f_s}{\partial x}\big(x_s^{\hat\varphi},\hat\varphi_{s}(x_s^{\hat\varphi})\big) + \frac{\partial f_s}{\partial u}\big(x_s^{\hat\varphi},\hat\varphi_s(x_s^{\hat\varphi})\big)\frac{\partial \hat\varphi_s}{\partial x}(x_s^{\hat\varphi})\Big)^\top \mathbb E\big[x_{T}^{{\hat\varphi}}-\mathfrak t(x_0)\big|x_{t-1}^{{\hat\varphi}}\big].
		\end{align}
		If $\hat\varphi$ is optimal for problem $(\mathcal P)$, then 
		\begin{enumerate}
			\item[(a)] ${p}_T =  x_T^{\hat\varphi} -\mathbb E\big[\mathfrak t(x_0)\,\big|\, x_{T-1}^{\hat\varphi}\big]$ and, for all $t\in\{1,\dots, T-1\}$,
			\begin{align}\label{MPXms}
				\mathbb E\big[{p}_t\,\big|\, x_{t}^{\hat\varphi}\big] = \frac{\partial f_t}{\partial x}\big(x_t^{\hat\varphi},\hat\varphi_t(x_t^{\hat\varphi})\big)^\top{p}_{t+1};
			\end{align}
			\item[(b)] for all $t\in\{0,\dots, T-1\}$,
			\begin{align}\label{MPYms}
				0\in \frac{\partial f_t}{\partial u}\big(x_t^{\hat\varphi},\hat\varphi_t(x_t^{\hat\varphi})\big)^\top {p}_{t+1}+ N_{U}\big(\hat\varphi_t(x_{t}^{\hat\varphi})\big).
			\end{align}
		\end{enumerate}
	\end{theorem}
	The variational inequality resulting from condition $(b)$ in Theorem \ref{Pontryagin} simplifies significantly in the absence of control constraints.
	\begin{corollary}
		Let the assumptions of Theorem \ref{Pontryagin} hold and suppose that  \( U = \mathbb{R}^m \). If \( \hat\varphi \in \Phi \) is optimal for problem \( (\mathcal{P}) \), then the sequence \( \{p_t\}_{t=1}^T \) in (\ref{adjoint}) satisfies
		\[
		\frac{\partial f_t}{\partial u}(x_t^{\hat\varphi}, \hat\varphi_t(x_t^{\hat\varphi}))^\top p_{t+1} = 0\quad \forall t \in \{0, \dots, T-1\}.
		\]
	\end{corollary}

	\subsection{A synthetic gradient for  descent methods}
	The variational principle in the last section suggest that is relevant to associate each Markov policy $\varphi \in \Phi$ with a sequence $\{p_t^{\varphi}\}_{t=1}^{T-1}$ of adjoint states and consider  
	\[
	\left\{\frac{\partial f_t}{\partial u} \big( x_t^{\varphi}, \varphi_t(x_t^{\varphi}) \big)^\top p_{t+1}^{\varphi} \right\}_{t=0}^{T-1}
	\]  
	as a representation of the gradient of the objective functional at $\varphi$. However, this expression is not well-suited for practical computation, particularly for gradient descent methods.  The main difficulty arises from the fact that Markov policies are sequences of functions rather than random variables. To address this issue, we define, for each $\varphi\in\Phi$, the sequence $\{\lambda_{t}^\varphi\}_{t=1}^{T}$ of functions $\lambda_t^\varphi:\mathbb R^n\to\mathbb R^n$ given by 
	\begin{align*}
		\lambda_t^\varphi(x):=\mathbb E\Big[\prod_{s=t}^{T-1}\Big(\frac{\partial f_s}{\partial x}\big(x_s^{\varphi},\varphi_s(x_s^{\varphi})\big) + \frac{\partial f_s}{\partial u}\big(x_s^{\varphi},\varphi_s(x_s^{\varphi})\big)\frac{\partial \varphi_s}{\partial x}(x_s^{\varphi})\Big)^\top
		\big(x_T^\varphi-\mathfrak t(x_0)\big)\,\Big|\, x_{t-1}^\varphi=x\Big].
	\end{align*}
	It is straightforward to verify that, for any  $\varphi\in\Phi$, the adjoint states satisfy
	\begin{align*}
		p_{t}^\varphi = \lambda_t^\varphi (x_{t}^{\varphi})\quad \forall t\in \{1,\dots, T\}.
	\end{align*}
	Now, we can define 	the \textit{synthetic gradient}  $\nabla \mathcal J(\varphi):=\big\{\big[\nabla \mathcal J(\varphi)\big]_t\big\}_{t=0}^{T-1}$ at Markov policy $\varphi\in\Phi$  as the sequence of functions $	\big[\nabla \mathcal J(\varphi)\big]_t:\mathbb R^n\to\mathbb R^m$ given by 
	\begin{align}\label{sgf}
		\big[\nabla \mathcal J(\varphi)\big]_t:= \frac{\partial f_t}{\partial u}\big(x,\varphi_t(x)\big)^\top \lambda_{t+1}^{\varphi}(x). 
	\end{align}
	With a well-defined notion of a gradient, we can formally devise a projected gradient descent algorithm for problem $(\mathcal{P})$ based on the iteration  
	\[
	\varphi^{i+1} = \Pi_{\Phi} \left( \varphi^i - \alpha \nabla \mathcal{J}(\varphi^i) \right),
	\]  
	where $\varphi^0 \in \Phi$ is a given initial policy, $\alpha > 0$ is the step size, and $\Pi_{\Phi}$ denotes the projection onto the feasible set $\Phi$, to be understood in the pointwise sense.
	
	\begin{algorithm}
		\caption{Projected synthetic gradient descent}
		\begin{algorithmic}[1]
			\Require Initial policy $\varphi^0 \in \Phi$, step size $\alpha > 0$, number of iterations $K$
			\State Initialize $\varphi \gets \varphi^0$
			
			\For{$i \gets 0$ to $K-1$}
			\State Compute gradient $\nabla \mathcal{J}(\varphi)$
			\State Update policy: ${\varphi} \gets \varphi - \alpha \nabla \mathcal{J}(\varphi)$
			\State Project onto feasible set: $\varphi \gets \Pi_{\Phi}({\varphi})$
			\EndFor
			
			\State \Return $\varphi$
		\end{algorithmic}
	\end{algorithm}


	\subsection{A couple of examples}
	We give a couple of simple examples to illustrate how the \textit{synthetic gradient} method described above can be used. The first focuses on a problem  where the policy updates can be computed analytically, and the second one on a problem where policy updates are computed numerically. 
	\subsubsection{Transporting a Gaussian}
	For illustration purposes, we consider a one-step distribution steering problem, which allows us to visualize the evolving Markov policies produced by the \textit{synthetic gradient} method. This setup is inspired by physical systems such as cold atoms in optical lattices or electrons in periodic potentials, where particles are subject to nonlinear forces. 
	\smallbreak 
	The state at time \(t\in\{0,1\}\) is denoted
	\[
	x_t = \begin{bmatrix} p_t \\ q_t \end{bmatrix} \in \mathbb{R}^2,
	\]
	where \(p_t\) and \(q_t\) are spatial coordinates. The system evolves according to
	\[
	x_1 = x_0 + f(x_0) + u \cdot G(x_0),
	\]
	where
	\[
	f(x_0) = \begin{bmatrix} \sin(p_0) \\ -\sin(q_0) \end{bmatrix}, \quad
	G(x_0) = \begin{bmatrix} 1 \\ 1 + \cos(p_0) \end{bmatrix}.
	\]
	The initial state is sampled from a Gaussian, i.e., $x_0\sim\mathcal{N}(m, I)$, where $m \in\mathbb R^d 
	$.
	The target distribution is \(\mathcal{N}(0, I)\). The control input \(u \in \mathbb{R}\) is applied once and determined via a policy \(u = \varphi(x_0)\), where   a Borel function \(\varphi: \mathbb{R}^2 \to \mathbb{R}\) is to be found. 
	\smallbreak 
	We define the transport map \(\mathfrak t(x) = x - m\), which applied to \(x_0 \sim \mathcal{N}(m, I)\) yields \(\mathfrak t(x_0)\sim\mathcal{N}(0, I)\). Let $\Phi$ denote the set of Borel functions $\varphi:\mathbb R^2\to\mathbb R$. The optimal control problem $(\mathcal P)$ then becomes 
	\[
	\min_{\varphi\in\Phi}\frac{1}{2} \, \mathbb{E} \left[ \left| x_1 - (x_0 - m) \right|^2 \right].
	\]
	This toy model captures key aspects of real quantum control problems, e.g., periodic potentials, nonlinear drift, and limited control.  It reflects a physical scenario where particles in a periodic field are initially displaced and a single, spatially varying control pulse is used to steer the ensemble toward a centered state.
	\smallbreak 
	Since the time horizon is $T=1$, the synthetic gradient simply becomes
	\[
	\left[ \nabla \mathcal J(\varphi) \right](x) = G(x)^\top \left( f(x) + \varphi(x) \cdot G(x) + m \right).
	\]
	\begin{figure}[h]
		\centering
		\includegraphics[width=0.8\textwidth]{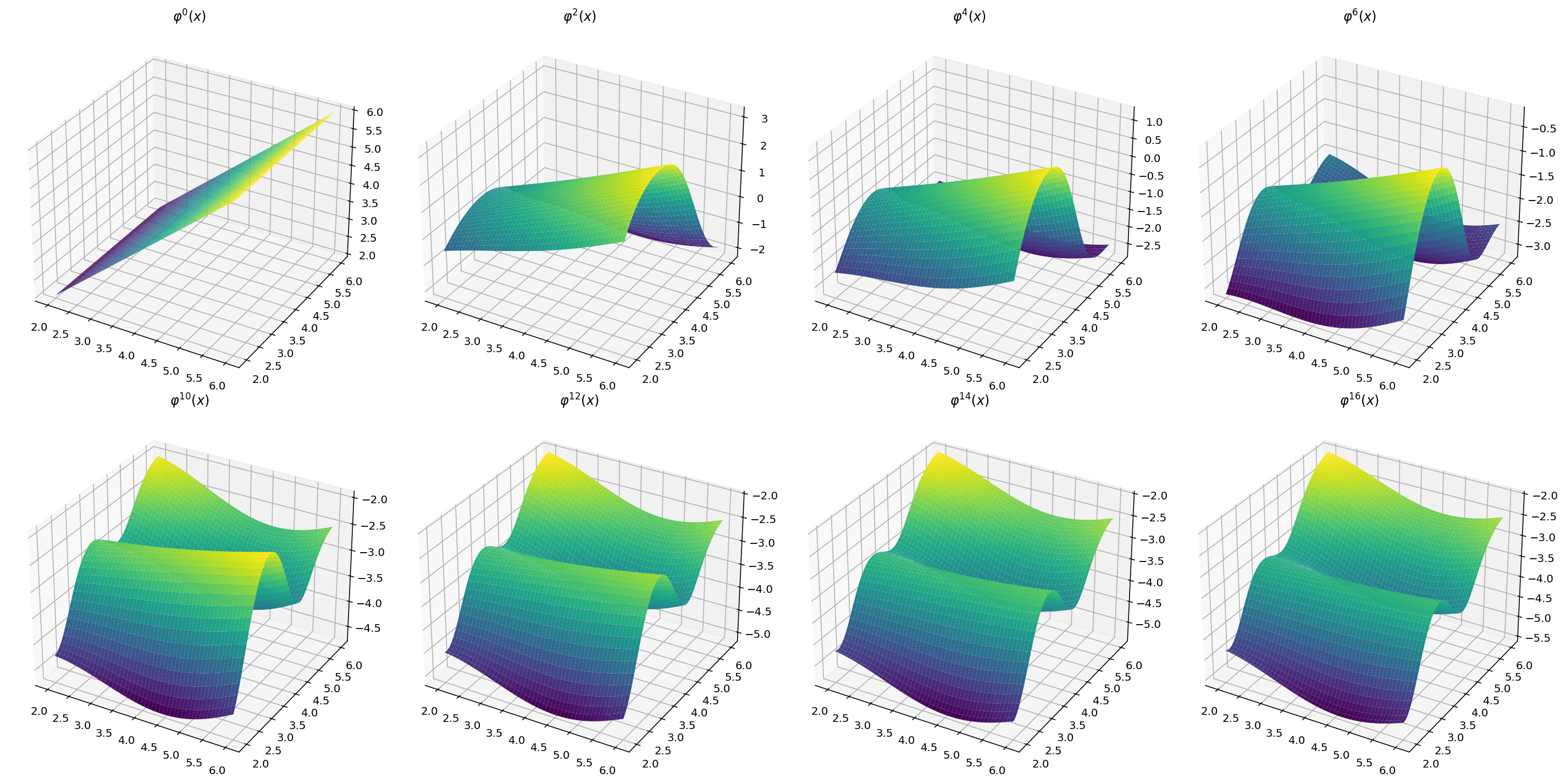}
		\caption{Evolution of policies}
		\label{fig:fig1}
	\end{figure}
	For a step-size $\alpha>0$, the update rule is then
	\[
	\varphi^{i+1}(x) = \varphi^i(x) - \alpha \cdot G(x)^\top \left( f(x) + \varphi^i(x) \cdot G(x) + m \right).
	\]
	This is a linear recursion in \(\varphi^i(x)\), and its closed-form solution is given by
	\[
	\varphi^i(x) 
	= \left(1 - \alpha |G(x)|^2\right)^i \cdot \varphi^0(x)
	- \frac{G(x)^\top (f(x) + m)}{|G(x)|^2} \left( 1 - \left(1 - \alpha |G(x)|^2\right)^i \right).
	\]
	If we take $\alpha < 2/5$,  then it can be seen that  $\left| 1 - \alpha |G(x)|^2 \right| < 1$ 
	for all $x\in\mathbb R^2$, which guarantees the pointwise convergence of the iteration. The pointwise limit is given by 
	\begin{align*}
		\hat \varphi(x) = -\frac{G(x)^\top (f(x) + \mu)}{|G(x)|^2}.
	\end{align*}
	Taking \(m = \left[\begin{smallmatrix} 4 \\ 4 \end{smallmatrix}\right]\), $\varphi^0(p,q)=p$ and \(\alpha = 0.15\), we plot in Figure~\ref{fig:fig1} the evolution of the policies \(\{\varphi^i\}_{i \in \mathbb{N}}\).  
	In Figure~\ref{fig:fig2}, we compare the distribution of \(x_1\) under the optimal policy \(\hat\varphi\) with the target one. Although the target distribution is \(\mathcal{N}(0, I)\), the resulting distribution does not  exactly match it due to the nonlinearities in the system.
	\begin{figure}[h]
		\centering
		\includegraphics[width=0.65\textwidth]{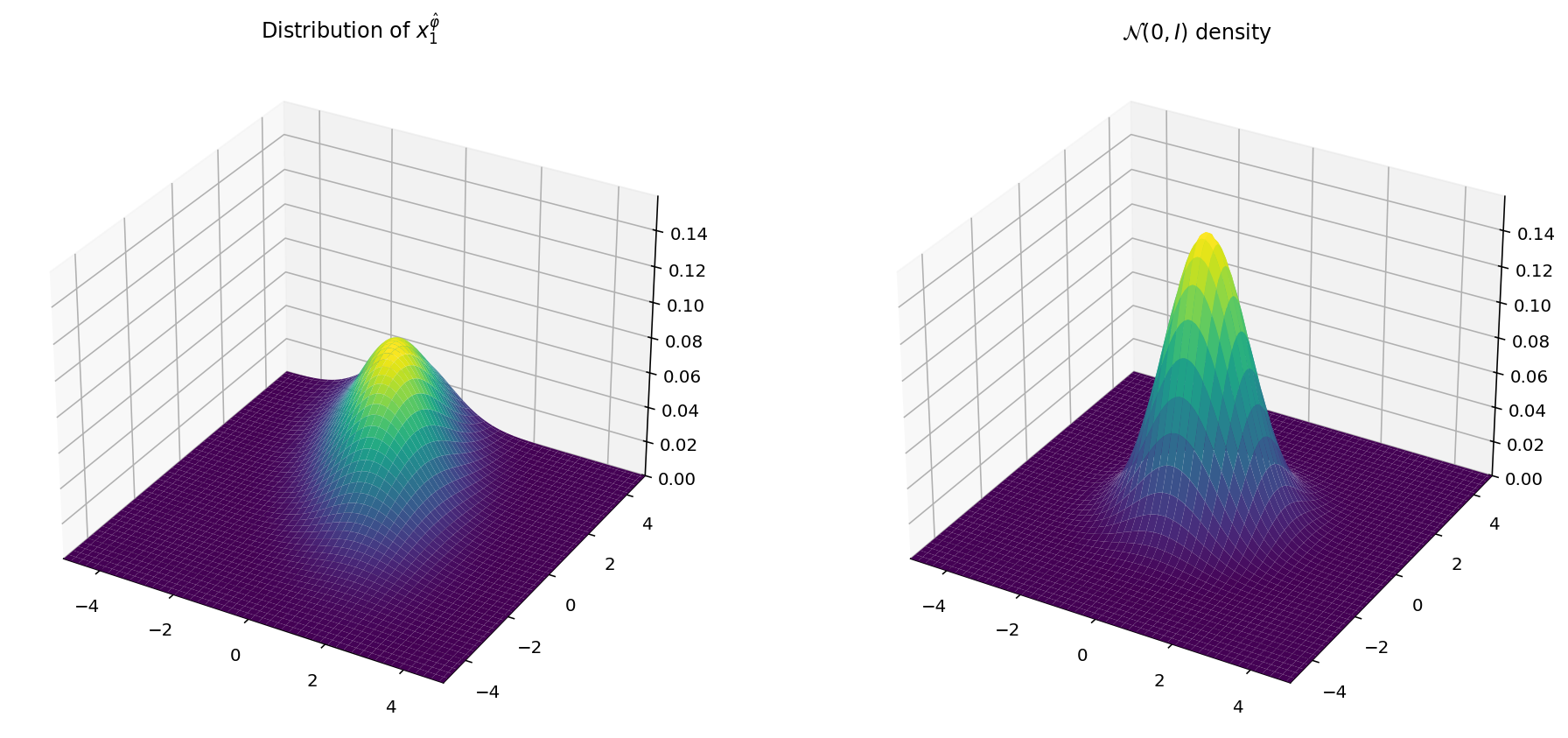}
		\caption{Comparison of distributions}
		\label{fig:fig2}
	\end{figure}
	
	\subsubsection{Collapsing a Gaussian}\label{Cog}
	In this example, we illustrate how the proposed gradient based method can be used to transport a Gaussian initial distribution toward the origin through a sequence of control actions.
	\begin{figure}[h]
		\centering
		\includegraphics[width=0.8\textwidth]{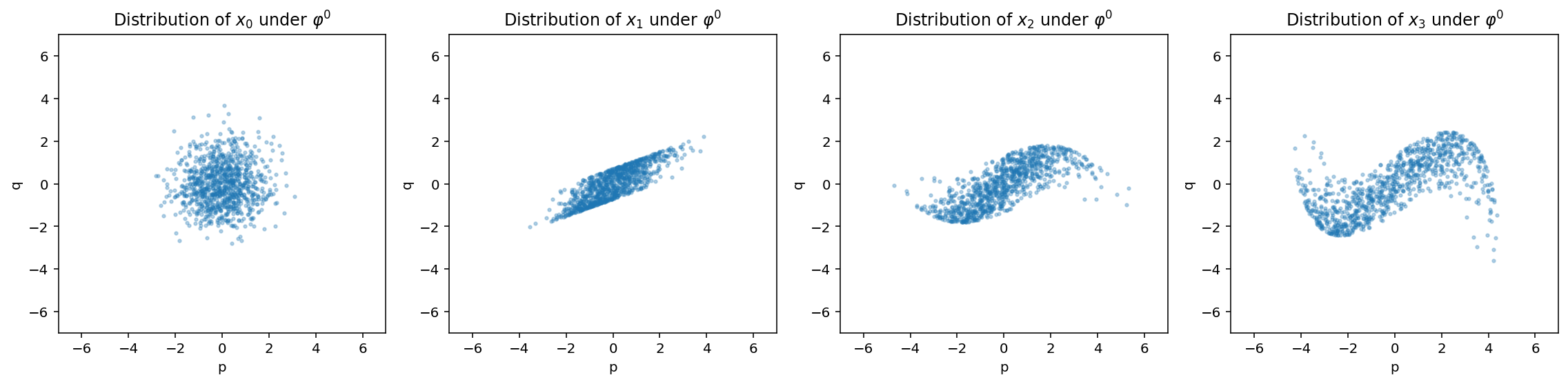}
		\caption{Samples of state distributions under the initial policy $\varphi^0$.}
		\label{fig:phi0}
	\end{figure}
	Let the state at time $t \in \{0,1,2,3\}$ be denoted by
	\[
	x_t = \begin{bmatrix} p_t \\ q_t \end{bmatrix} \in \mathbb{R}^2.
	\]
	The system evolves according to the nonlinear dynamics
	\[
	f(x, u) = \begin{bmatrix} p + q + u_1 \\ \beta \cdot q + \sin(p) + u_2 \end{bmatrix}, \quad u = \begin{bmatrix} u_1 \\ u_2 \end{bmatrix} \in \mathbb{R}^2,
	\]
	where $f : \mathbb{R}^2 \times \mathbb{R}^2 \to \mathbb{R}^2$ is the transition function and $\beta = 0.9$ is a damping parameter. The system is fully actuated, the control enters linearly in each coordinate.
	\smallbreak 
	The transport map is chosen as $\mathfrak{t}(x) = 0$. In this setting, the optimal control problem becomes
	\[
	\min_{\varphi\in\Phi} \frac{1}{2} \mathbb{E} \left[ | x_T |^2 \right].
	\]
	The aim is to minimize this quantity, with $T=3$,  by iteratively updating a time-dependent closed-loop policies $\varphi_0,\varphi_1,\varphi_2 : \mathbb{R}^2 \to \mathbb{R}^2$ $\varphi : \mathbb{R}^2 \to (\mathbb{R}^2)^T$,  using the {gradient-based} method described in the previous section.
	We initialize the policy with $\varphi_t^0(x) := -0.5x$ for each $t = 0, 1, 2$. Updates of the form
	\[
	\varphi^{(k+1)}(x) = \varphi^{(k)}(x) - \alpha \cdot \nabla \mathcal J(\varphi^{(k)})(x),
	\]
	are carried out, where $\alpha = 0.14$ is a fixed step size and $K=3$. The synthetic gradient $\nabla \mathcal J(\varphi^{(k)})(x)$ is computed with formula (\ref{sgf}) , using automatic differentiation to evaluate the Jacobians of the current policy $\varphi^{(k)}$. 
	
\begin{figure}[h]
		\centering
		\includegraphics[width=0.8\textwidth]{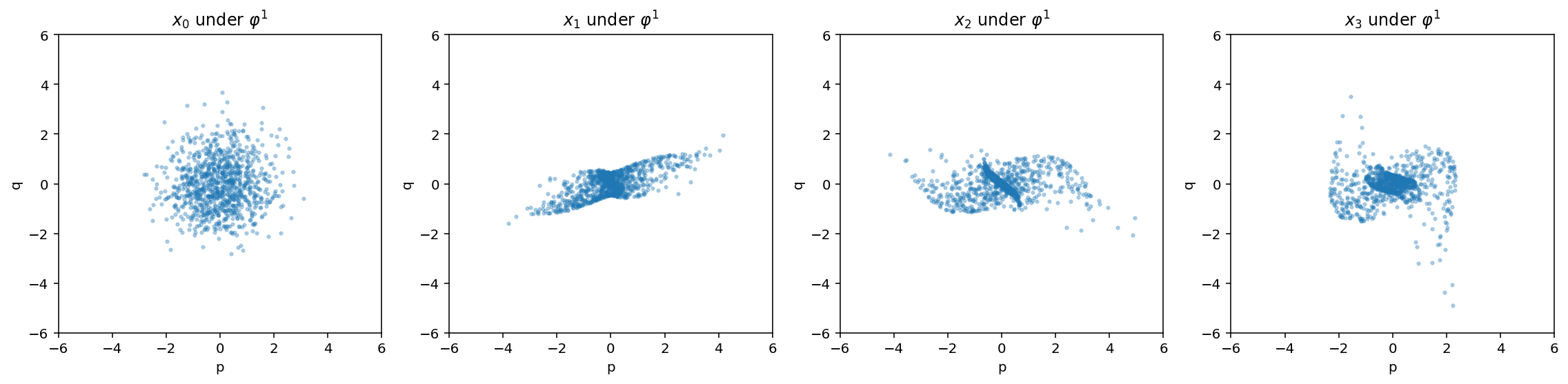}
		\caption{Samples of state distributions under the updated policy $\varphi^1$.}
		\label{fig:phi1}
	\end{figure}

\begin{figure}[h]
		\centering
		\includegraphics[width=0.8\textwidth]{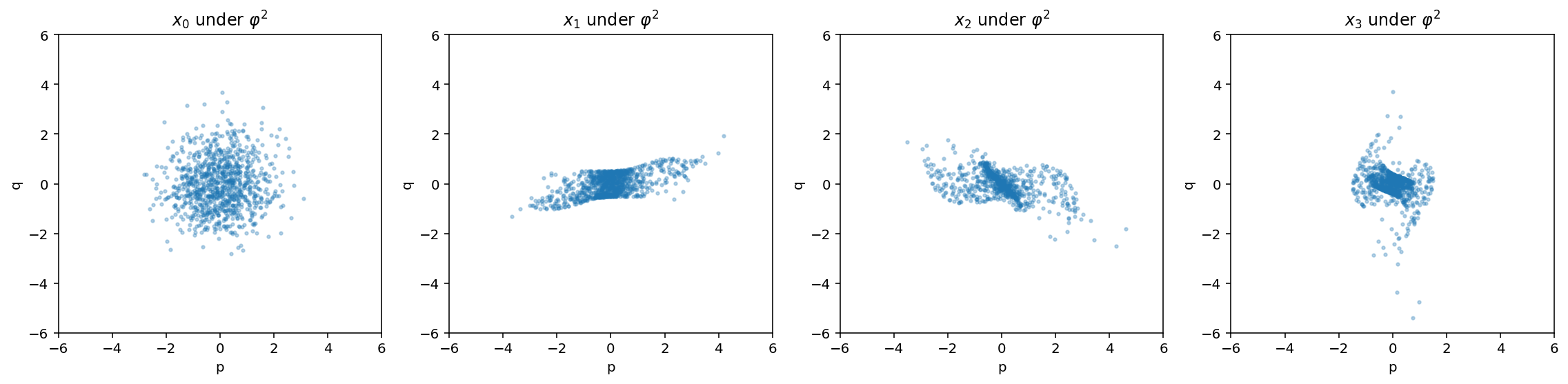}
		\caption{Samples of state distributions under the updated policy $\varphi^2$.}
		\label{fig:phi2}
	\end{figure}

\begin{figure}[h]
		\centering
		\includegraphics[width=0.8\textwidth]{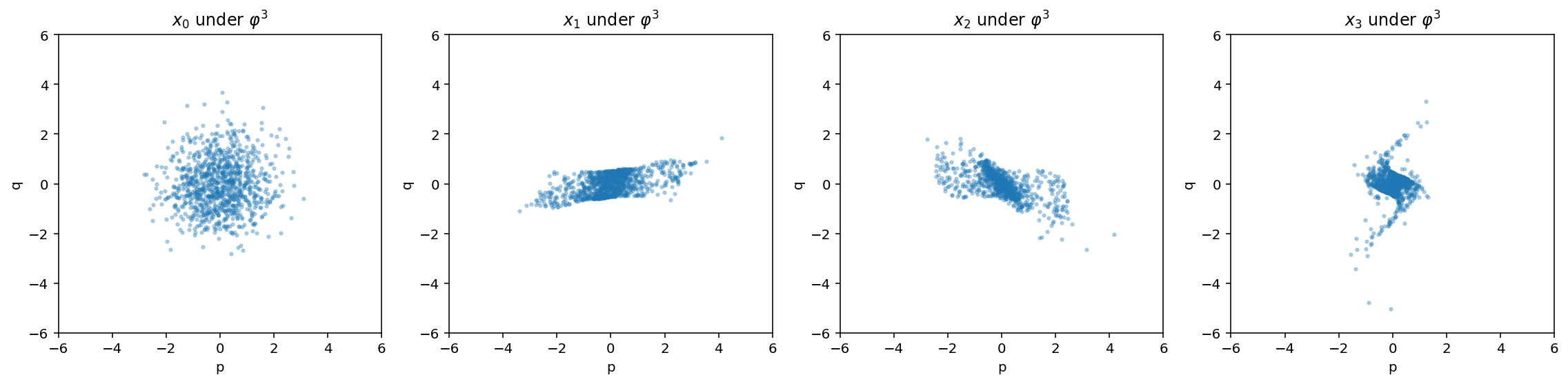}
		\caption{Samples of state distributions under the updated policy $\varphi^3$.}
		\label{fig:phi3}
	\end{figure}

	\smallbreak 
	We sample $x_0 \sim \mathcal{N}(0, I)$ to carry out numerical computations. Figures~\ref{fig:phi0}--\ref{fig:phi3} show the empirical distributions of the states $x_t$ at each time step $t = 0,1,2,3$ under successive policies $\varphi^k$.
	
	To summarize the numerical results, Table~\ref{tab:objective-values} reports the value of the objective functional $\mathcal{J}(\varphi^k)$ across iterations. We observe a rapid decrease in cost after the first update, followed by smaller refinements in subsequent steps. 
	
	\begin{table}[h]
		\centering
		\renewcommand{\arraystretch}{1.2}
		\begin{tabular}{c|c}
			Iteration $k$ & Objective $\mathcal{J}(\varphi^k)$ \\ \hline
			0 & 3.1826 \\
			1 & 0.7327 \\
			2 & 0.3506 \\
			3 & 0.2813 \\
		\end{tabular}
		\caption{Values of the objective $\mathcal{J}(\varphi^k)$ over successive iterations.}
		\label{tab:objective-values}
	\end{table}


	\section{Proofs}
	\subsection{Proof of Proposition \ref{Propp}}
	Due to the hypotheses, we can employ a standard selection theorem, e.g., \cite[Proposition D.6-(a)]{HL_1996}, to conclude that, for each $t\in\{0,\dots, T-1\}$, there exists a Borel measurable function $\varphi_t:\mathbb R^n\to\mathbb R^m$ such that 
	\begin{align*}
		\inf_{u\in U} |f_t(x,u)| = |f_t(x,\varphi_t(x))| \quad \forall x\in\mathbb R^d.  
	\end{align*}
	Set $\varphi:=\{\varphi_t\}_{t=0}^{T-1}$; we will prove that $\mathcal J(\varphi)<+\infty$. 
	\smallbreak 
	We prove by mathematical  induction that $x_t^\varphi$ has finite second moment for every $t\in\{0,\dots, T\}$. By assumption, $x_0$ has  finite  second moment. Now, let $\tau\in\{1\dots, T-1\}$ and suppose that $x_\tau^\varphi$ has finite second moment. Then, due to hypothesis $(ii)$,
	\begin{align*}
		\mathbb E|x_{\tau+1}^\varphi|^2 = \mathbb E|f_{t}\big(x_{\tau}^\varphi, \varphi_\tau(x_{\tau}^\varphi)\big)|^2 \le m_\tau^2\mathbb E \big(1+|x_{\tau}^\varphi|\big)^2<+\infty. 
	\end{align*}
	We conclude that $x_{\tau+1}^\varphi$ has finite second moment;  thereby completing the induction argument.
	\smallbreak 
	Since both  $x_{T}^\varphi$ and $\mathfrak t(x_0)$ have a finite second moment, we conclude that\begin{align*}
		\mathcal J(\varphi) = \frac{1}{2} \mathbb E\big | x_T^{\varphi} - \mathfrak t(x_0)\big|^2 <+\infty. 
	\end{align*}

	\subsection{Proof of Theorem \ref{Pontryagin}}
	For each $t\in\{1,\dots, T-1\}$, consider the function $q_t:\mathbb R^n\to \mathbb R^n$ given by 
	\begin{align*}
		q_t(x):=\frac{\partial f_t}{\partial x}\big(x,\hat\varphi_t(x)\big) + \frac{\partial f_t}{\partial u}\big(x,\hat\varphi_t(x)\big)\frac{\partial \hat\varphi_t}{\partial x}(x).
	\end{align*}
	Now, set  ${p}_T:= \mathbb E\big[x_{T}^{{\hat\varphi}}-\mathfrak t(x_0)\big|x_{T-1}^{{\hat\varphi}}\big]$ and, for $t\in\{1,\dots, T-1\}$,
	\begin{align*}
		{p}_{t}:= \Big[\prod_{s=t}^{T-1} q_s(x_s^{\hat\varphi})\Big]^\top \mathbb E\big[x_{T}^{{\hat\varphi}}-\mathfrak t(x_0)\big|x_{t-1}^{{\hat\varphi}}\big].
	\end{align*}	
	We prove below that $\{{p}_t\}_{t=0}^{T-1}$ satisfies the desired properties. Let $\tau\in \{0,\dots, T-1\}$ be given. 
	\smallbreak
	\noindent
	\textbf{Step 1} \textit{(An auxiliary dynamical system)}.
	For any given Borel function $\psi:\mathbb R^n\to U$, we consider the sequence $\{z_{t}^\psi\}_{t=\tau+1}^{T}$  given by 
    \begin{align*}
        z_{\tau+1}^\psi=\frac{\partial f_\tau}{\partial u}\big(x_{\tau}^{\hat\varphi},\hat\varphi_\tau(x_{\tau}^{\hat\varphi})\big)\big(\psi(x_{\tau}^{\hat\varphi})-\hat\varphi(x_{\tau}^{\hat\varphi})\big)\quad\text{and}\quad z_{t+1}^{\psi} = q_{t}(x_{t}^{\hat\varphi})z_{t}^\psi.
    \end{align*}
	We can easily see that, for every $t\in\{\tau+1,\dots, T\}$, 
	\begin{align*}
		z_{t}^\psi = \Big[\prod_{s=\tau+1}^{t-1}q_{s}(x_s^{\hat\varphi}) \Big]\frac{\partial f_\tau}{\partial u}\big(x_{\tau}^{\hat\varphi},\hat\varphi_\tau(x_{\tau}^{\hat\varphi})\big)\big(\psi(x_{\tau}^{\hat\varphi})-\hat\varphi(x_{\tau}^{\hat\varphi})\big).
	\end{align*}
	\textbf{Step 2} \textit{(Calculation of a Gateaux differential)}.  Given a Borel function $\psi:\mathbb R^n\to U$, consider the sequence $\varphi^{\tau,\psi}:=\{\varphi^{\tau,\psi}_t\}_{t=0}^{T-1}$ given by 
	\begin{align}
		\varphi_t^{\tau,\psi}(x):=\left\{\begin{array}{lcc}
			\psi(x) - \hat\varphi_\tau(x) & \text{if} & t=\tau\\
			0 & \text{if} & t\neq \tau.
		\end{array}
		\right. 
	\end{align}
	Due to convexity of the control set $U$, ${\hat\varphi} + \varepsilon\varphi^{\tau,u}\in \Phi$ for all $\varepsilon\in[0,1]$. We prove below that
	\begin{align}\label{gade}
		\lim_{\varepsilon\longrightarrow0^+}	\frac{\mathcal J({\hat\varphi} + \varepsilon\varphi^{\tau,\psi}) - \mathcal J({\hat\varphi})}{\varepsilon} = \mathbb E\Big\langle  \frac{\partial f_{\tau}}{\partial u}\big(x_\tau^{{\hat\varphi}}, \hat\varphi_\tau(x_{\tau}^{{\hat\varphi}})\big) ^\top{p}_{\tau+1},\psi(x_{\tau}^{\hat\varphi}) - {\hat\varphi_\tau}(x_{\tau}^{{\hat\varphi}})\Big\rangle. 
	\end{align}
	Due the hypotheses, each $q_{t}$ is uniformly continuous, and hence we can find a continuous nondecreasing function  $\omega:[0,1]\to\mathbb R$ vanishing at zero such that
	\begin{align*}
		| x_{T}^{\hat\varphi + \varepsilon\varphi^{\tau,\psi}} - x_{T}^{\hat\varphi} - \varepsilon z_{T}^\psi|\le \varepsilon \omega(\varepsilon)|\hat\varphi_\tau(x_{\tau}^{\hat\varphi}) - \psi(x_{\tau}^{\hat\varphi})|\quad \forall\varepsilon\in[0,1]. 
	\end{align*}
	The previous estimate allows to interchange the derivative and expectation signs. Then, 
	\begin{align*}
		&\frac{d}{d\varepsilon}\mathcal J\big({\hat\varphi}+\varepsilon\varphi^{\tau,\psi}\big)\Big|_{\varepsilon=0} = \mathbb E\big\langle x_{T}^{\hat\varphi} - \mathfrak t(x_0), z_{T}^{\psi}\big\rangle\\
		&\,\,\,=\mathbb E\Big\langle x_{T}^{\hat\varphi} - \mathfrak t(x_0), \Big[\prod_{s=\tau+1}^{T-1}q_{s}(x_s^{\hat\varphi}) \Big]\frac{\partial f_\tau}{\partial u}\big(x_{\tau}^{\hat\varphi},\hat\varphi_\tau(x_{\tau}^{\hat\varphi})\big)\big(\psi(x_{\tau}^{\hat\varphi})-\hat\varphi_\tau(x_{\tau}^{\hat\varphi})\big)\Big\rangle\\
		&\,\,\,= \mathbb E\Big\langle \mathbb E\big[x_{T}^{\hat\varphi} - \mathfrak t(x_0)\,\big|\, x_{\tau}^{\hat\varphi}\big], \Big[\prod_{s=\tau+1}^{T-1}q_{s}(x_s^{\hat\varphi}) \Big]\frac{\partial f_\tau}{\partial u}\big(x_{\tau}^{\hat\varphi},\hat\varphi_\tau(x_{\tau}^{\hat\varphi})\big)\big(\psi(x_{\tau}^{\hat\varphi})-\hat\varphi_\tau(x_{\tau}^{\hat\varphi})\big)\Big\rangle
	\end{align*}
	Whence (\ref{gade}) follows directly.
	\smallbreak 
	\textbf{Step 3} \textit{(First-order necessary condition)}.  By the previous step, for any $\psi:\mathbb R^n\to U$ Borel, 
	\begin{align*}
		0&\le 	\lim_{\varepsilon\longrightarrow0^+}	\frac{\mathcal J({\hat\varphi} + \varepsilon\varphi^{\tau,{u}}) - \mathcal J({\hat\varphi})}{\varepsilon} \\
        &= \mathbb E\Big \langle \frac{\partial f_\tau}{\partial u}\big(x_{\tau}^{{\hat\varphi}},\hat\varphi_\tau(x_{\tau}^{{\hat\varphi}})\big)^\top{p}_{\tau+1}, \psi(x_{\tau}^{\hat\varphi}) - \hat\varphi_\tau(x_{\tau}^{{\hat\varphi}}) \Big\rangle.
	\end{align*}
	We can then conclude that, for any $u\in U$, 
	\begin{align}\label{FONCe}
		\Big \langle \frac{\partial f_\tau}{\partial u}\big(x_{\tau}^{{\hat\varphi}},\hat\varphi_\tau(x_{\tau}^{{\hat\varphi}})\big)^\top{p}_{\tau+1}, u - \hat\varphi_\tau(x_{\tau}^{{\hat\varphi}}) \Big\rangle\ge0.
	\end{align}
	By definition of normal cone, we can write this as 
	\begin{align*}
		0\in \frac{\partial f_\tau}{\partial u}\big(x_{\tau}^{{\hat\varphi}},\hat\varphi_\tau(x_{\tau}^{{\hat\varphi}})\big)^\top {p}_{\tau+1} + N_{U}\big({\hat\varphi}(x_{\tau}^{{\hat\varphi}})\big).
	\end{align*}
	\textbf{Step 4} \textit{(The recurrence of the adjoint)}. Observe that $\hat\varphi_\tau(x_{\tau}^{{\hat\varphi}}+\varepsilon v)\in U$ for any $\varepsilon>0$ and  $v\in\mathbb R^d$.  Combining this with (\ref{FONCe}), for any $v\in\mathbb R^n$,
	\begin{align*}
		\Big \langle \frac{\partial f_\tau}{\partial u}\big(x_{\tau}^{{\hat\varphi}},\hat\varphi_\tau(x_{\tau}^{{\hat\varphi}})\big)^\top{p}_{\tau+1}, \frac{\hat\varphi_\tau(x_{\tau}^{{\hat\varphi}}+\varepsilon v)- \hat\varphi_\tau(x_{\tau}^{{\hat\varphi}})}{\varepsilon} \Big\rangle\ge0 \quad \forall \varepsilon>0. 
	\end{align*}
	Taking limit as $\varepsilon\longrightarrow0^+$, we obtain, for any $v\in\mathbb R^n$, 
	\begin{align*}
		\Big \langle \frac{\partial f_\tau}{\partial u}\big(x_{\tau}^{{\hat\varphi}},\hat\varphi_\tau(x_{\tau}^{{\hat\varphi}})\big)^\top{p}_{\tau+1}, \frac{\partial \hat\varphi_\tau}{\partial x}(x_{\tau}^{{\hat\varphi}})v \Big\rangle\ge0.
	\end{align*}
	Since this holds for $v\in\mathbb R^d$ arbitrary, we conclude that 
	\begin{align}\label{Orthadjo}
		\frac{\partial \hat\varphi_\tau}{\partial x}(x_{\tau}^{{\hat\varphi}})^\top  \frac{\partial f_\tau}{\partial u}\big(x_{\tau}^{{\hat\varphi}},\hat\varphi_\tau(x_{\tau}^{{\hat\varphi}})\big)^\top {p}_{\tau +1} = 0. 
	\end{align}
	Finally, employing (\ref{Orthadjo}), we see that
	\begin{align*}
		\mathbb E\big[{p}_\tau\,\big|\,x_{\tau}^{\hat\varphi}\big] &= \mathbb E \Big[\Big[\prod_{s=\tau}^{T-1} q_s(x_s^{\hat\varphi})\Big]^\top \mathbb E\big[x_{T}^{{\hat\varphi}}-\mathfrak t(x_0)\big|x_{\tau-1}^{{\hat\varphi}}\big]\,\Big|\, x_{\tau}^{\hat{\varphi}}\Big]\\
        &=  \Big[\prod_{s=\tau}^{T-1} q_s(x_s^{\hat\varphi})\Big]^\top \mathbb E\big[x_{T}^{{\hat\varphi}}-\mathfrak t(x_0)\big|x_{\tau}^{{\hat\varphi}}\big] \\
		&=q_\tau(x_\tau^{\hat\varphi})^\top\Big[\prod_{s=\tau+1}^{T-1} q_s(x_s^{\hat\varphi})\Big]^\top \mathbb E\big[x_{T}^{{\hat\varphi}}-\mathfrak t(x_0)\big|x_{\tau}^{{\hat\varphi}}\big]=q_\tau(x_\tau^{\hat\varphi})^\top p_{\tau+1}\\
		&=\Big(\frac{\partial f_\tau}{\partial x}\big(x_{\tau}^{{\hat\varphi}},\varphi(x_{\tau}^{{\hat\varphi}})\big)^\top + \frac{\partial \hat\varphi_\tau}{\partial x}(x_{\tau}^{{\hat\varphi}})^\top\frac{\partial f_\tau}{\partial u}\big(x_{\tau}^{{\hat\varphi}},\varphi(x_{\tau}^{{\hat\varphi}})\big)^\top\Big)  p_{\tau+1}\\
        &=\frac{\partial f_\tau}{\partial x}\big(x_{\tau}^{{\hat\varphi}},\varphi(x_{\tau}^{{\hat\varphi}})\big)^\top p_{\tau+1}.
	\end{align*}
	This completes the proof. 

\begin{thebibliography}{99}
		
		\bibitem{AM_1971} B. D. O. Anderson and J. B. Moore, {\it Linear optimal control}, Prentice-Hall, Englewood Cliffs, NJ, 1971.
		
		\bibitem{AlNassirLafta_2021}
		S.~Al-Nassir and A.~Lafta.
		\newblock Analysis of a harvested discrete-time biological model.
		\newblock {\em Int. J. Nonlinear Anal. Appl.}, 12(2):2235--2246, 2021.
		
		
		\bibitem{Aseev_2017}
		S. M. Aseev, M. I. Krastanov and V. M. Veliov, Optimality conditions for discrete-time optimal control on infinite horizon, Pure Appl. Funct. Anal. {\bf 2} (2017), no.~3, 395--409.
		
		\bibitem{AL_1983} S. Aubry and P. Y. Le Daeron, {\it The discrete Frenkel-Kontorova model and its extensions I}, Physica D {\bf 8} (1983), 381--422.
		
		\bibitem{Bakolas_2016}
		E. Bakolas, Optimal covariance control for discrete-time stochastic linear systems subject to constraints, Proc. 55th IEEE Conf. on Decision and Control (CDC), Las Vegas, NV, USA, 2016, pp. 1153--1158.
		
		\bibitem{Bakolas_2017}
		E. Bakolas, Covariance control for discrete-time stochastic linear systems with incomplete state information, Proc. American Control Conference (ACC), Seattle, WA, USA, 2017, pp. 432--437.
		
		\bibitem{Bakolas_SCL_2018}
		E. Bakolas, Constrained minimum variance control for discrete-time stochastic linear systems, Systems Control Lett. {\bf 113} (2018), 109--116.
		
		\bibitem{Bakolas_Automatica_2018}
		E. Bakolas, Finite-horizon covariance control for discrete-time stochastic linear systems subject to input constraints, Automatica {\bf 91} (2018), 61--68.
		
		\bibitem{BalciBakolas_2021}
		I. M. Balci and E. Bakolas, Covariance control of discrete-time Gaussian linear systems using affine disturbance feedback control policies, Proc. 60th IEEE Conf. on Decision and Control (CDC), Austin, TX, USA, 2021, pp. 2324--2329.
		
		\bibitem{Bonnet_2019}
		B.~Bonnet.
		\newblock A Pontryagin maximum principle in Wasserstein spaces for constrained optimal control problems.
		\newblock {\em ESAIM Control Optim. Calc. Var.}, 25:52, 2019.
		
		
		
		\bibitem{BonnetFrankowska}
		B. Bonnet-Weill and H. Frankowska, Necessary optimality conditions for optimal control problems in Wasserstein spaces, Appl. Math. Optim. {\bf 84} (2021), S1281--S1330.
		
		\bibitem{BonnetRossi_2019}
		B.~Bonnet and F.~Rossi.
		\newblock The Pontryagin maximum principle in the Wasserstein space.
		\newblock {\em Calc. Var. Partial Differ. Equ.}, 58:11, 2019.
		
		\bibitem{BZ_2005} J. M. Borwein and Q. J. Zhu, {\it Techniques of variational analysis}, CMS Books in Mathematics, vol. 20, Springer-Verlag, New York, 2005.
		
		\bibitem{deBaydinEtAl}
		M. H. de Badyn et al., Discrete-Time Linear-Quadratic Regulation via Optimal Transport, 60th IEEE Conference on Decision and Control (CDC), Austin, TX, USA, 2021, pp. 3060-3065, doi: 10.1109/CDC45484.2021.9682825.
		
		\bibitem{ChenGeorgiouPavon_JOTA_2016}
		Y.~Chen, T.~T. Georgiou, and M.~Pavon.
		\newblock On the relation between optimal transport and Schrödinger bridges: A stochastic control viewpoint.
		\newblock {\em J. Optim. Theory Appl.}, 169(2):671--691, 2016.
		
		\bibitem{ChenGeorgiouPavon_TAC_PartI_2016}
		Y.~Chen, T.~T. Georgiou, and M.~Pavon.
		\newblock Optimal steering of a linear stochastic system to a final probability distribution, Part I.
		\newblock {\em IEEE Trans. Automat. Control}, 61(5):1158--1169, 2016.
		
		\bibitem{ChenGeorgiouPavon_TAC_PartII_2016}
		Y.~Chen, T.~T. Georgiou, and M.~Pavon.
		\newblock Optimal steering of a linear stochastic system to a final probability distribution, Part II.
		\newblock {\em IEEE Trans. Automat. Control}, 61(5):1170--1180, 2016.
		
		\bibitem{ChenGeorgiouPavon_TAC_PartIII_2018}
		Y.~Chen, T.~T. Georgiou, and M.~Pavon.
		\newblock Optimal steering of a linear stochastic system to a final probability distribution, Part III.
		\newblock {\em IEEE Trans. Automat. Control}, 63(9):3112--3118, 2018.
		
		
		\bibitem{CorellaHL_2019}
		A. Dom\'inguez Corella and O. Hern\'andez-Lerma, The maximum principle for discrete-time control systems and applications to dynamic games, J. Math. Anal. Appl. {\bf 475} (2019), no.~1, 253--277.
		
		
		\bibitem{HL_1996} O. Hern\'andez-Lerma and J.-B. Lasserre, {\it Discrete-time Markov control processes}, Applications of Mathematics (New York), 30, Springer, New York, (1996).
		
		\bibitem{ItoKashima_2025}
		K. Ito and K. Kashima, Maximum entropy density control of discrete-time linear systems with quadratic cost, IEEE Trans. Automat. Control {\bf 70} (2025), no.~5, 3024--3039.
		
		\bibitem{KeskinSaglam_2022}
		K.~Keskin and Ç.~Sağlam.
		\newblock Fatigue accumulation in dynamic contests.
		\newblock {\em Oper. Res. Lett.}, 50(3):268--273, 2022.

             \bibitem{KrastanovStefanov_2024}
M. I. Krastanov and B. K. Stefanov, A sufficient condition for a discrete-time optimal control problem, in {\it Large-scale scientific computations}, Lecture Notes in Comput. Sci., vol.~13952, Springer, Cham, 2024, pp.~184--192.

		
		\bibitem{P_2013} J.-P. Penot, {\it Calculus without derivatives}, Graduate Texts in Mathematics, vol. 266, Springer, New York, 2013.
		
		\bibitem{TerpinEtAl}
		A. Terpin, N. Lanzetti and F.~A. Dorfler, Dynamic programming in probability spaces via optimal transport, SIAM J. Control Optim. {\bf 62} (2024), no.~2, 1183--1206.

   
		
	\end{thebibliography}
\end{document}